\documentclass[11pt]{amsart}

\usepackage{amsmath}
\usepackage{amssymb}
\usepackage{mathrsfs}
\usepackage[all]{xy}
\usepackage{comment}
\usepackage{color}
\usepackage[colorlinks,citecolor=blue,urlcolor=black,linkcolor=black]{hyperref}

\newtheorem{theorem}{Theorem}[section]
\newtheorem{claim}[theorem]{Claim}

\newtheorem{lemma}[theorem]{Lemma}

\newtheorem{corollary}[theorem]{Corollary}

\theoremstyle{definition}
\newtheorem{definition}[theorem]{Definition}
\newtheorem{axiom}[theorem]{Axiom}

\newtheorem{question}[theorem]{Question}

\theoremstyle{remark}
\newtheorem{remark}[theorem]{Remark}

\newcount\skewfactor
\def\mathunderaccent#1#2 {\let\theaccent#1\skewfactor#2
\mathpalette\putaccentunder}
\def\putaccentunder#1#2{\oalign{$#1#2$\crcr\hidewidth
\vbox to.2ex{\hbox{$#1\skew\skewfactor\theaccent{}$}\vss}\hidewidth}}
\def\name{\mathunderaccent\tilde-3 }

\def\smallbox#1{\leavevmode\thinspace\hbox{\vrule\vtop{\vbox
   {\hrule\kern1pt\hbox{\vphantom{\tt/}\thinspace{\tt#1}\thinspace}}
   \kern1pt\hrule}\vrule}\thinspace}

\newcommand{\cf}{{\rm cf}}

\newcommand{\stick}{\ensuremath\mspace{2mu}\mid \mspace{-12mu}{\raise 0.4em \hbox{$\bullet$}}}

\def\qedref#1{$\qed_{\reforiginal{#1}}$}

\title{Echeloned saturation and forcing axioms}
\author{Shimon Garti}
\address{Institute of Mathematics,
 The Hebrew University of Jerusalem,
 Jerusalem 91904, Israel}
\email{shimon.garty@mail.huji.ac.il}
\thanks{}


\begin{document}
\let\labeloriginal\label
\let\reforiginal\ref

\begin{abstract}
We show that under the failure of Chang's conjecture and the assumption that Martin's axiom holds and the continuum is above $\omega_2$ there are no weakly Laver ideals over $\aleph_1$. This gives an answer to a question of Paul Larson. Likewise, we prove that under Baumgartner's axiom there is no Laver ideal over $\aleph_2$. This statement generalizes Larson's result at $\aleph_1$ under Martin's axiom.
\end{abstract}

\maketitle

\newpage

\section{Introduction}

Small accessible cardinals, especially $\aleph_1$ and $\aleph_2$, tend to behave like large cardinals under forcing axioms. This proclivity is manifested in several ways, some of which are related to saturated ideals. For example, it is shown in \cite{MR924672} that $\aleph_1$ carries a $\sigma$-complete $\aleph_2$-saturated ideal under Martin's maximum (ideals with these properties will be called Kunen ideals henceforth).
In fact, Martin's maximum causes ${\rm NS}_{\omega_1}$ to be a Kunen ideal.
In this light, the following theorem of Larson, \cite{MR2146222}, is a bit surprising:

\begin{theorem}
\label{thmlarson} Assume Martin's axiom and $2^\omega=\aleph_2$. \newline
Then there is no $\aleph_1$-complete $(\aleph_2,\aleph_2,\aleph_0)$-saturated ideal over $\aleph_1$.
\end{theorem}

\hfill \qedref{thmlarson}

The saturation property mentioned in this theorem (escorted by the additional property of $\aleph_1$-completeness) is called Laverness.
Laver, \cite{MR673792}, proved that if there is a huge cardinal then such an ideal can be forced over $\aleph_1$.
Nonetheless, Martin's axiom and $2^\omega=\omega_2$ (a fortiori, the proper forcing axiom or Martin's maximum) excludes Laver ideals from $\aleph_1$.
By the results of the current paper, Martin's axiom with $2^\omega>\omega_1$ has the same effect, provided that Chang's conjecture fails.

One moral of Larson's theorem is that Kunen's saturation and Laver's saturation differ significantly. A natural question raised by Larson in \cite{MR2146222} is whether $\aleph_1$ carries an $(\aleph_2,\aleph_1,\aleph_0)$-saturated ideal under some forcing axiom. Such an ideal (which will be called an echeloned ideal or weakly Laver) has a saturation property between Kunen's saturation and Laver's saturation. Larson asked, specifically, about the ideal of non-stationary sets ${\rm NS}_{\omega_1}$. The main result of this paper is the following:

\begin{theorem}
\label{thmthm} Assume Martin's axiom and $2^\omega>\omega_2$. Assume, further, that Chang's conjecture fails.
Then there are no weakly Laver ideals over $\aleph_1$.
\end{theorem}

\hfill \qedref{thmthm}

Unlike the assumption $2^\omega=\omega_2$ in Larson's theorem, here we assume that $2^\omega>\omega_2$.
It seems that a deep point lies here.
One aspect of this point is related to strong forcing axioms like \textsf{MM} or \textsf{PFA}.
Both imply Martin's axiom, and both imply $2^\omega=\omega_2$.\footnote{The fact that \textsf{MM} implies $2^\omega=\omega_2$ appears in \cite{MR924672}. It is also a consequence of \textsf{PFA} as proved by Veli\u{c}kovi\'c and Todor\v{c}evi\'c, see \cite{MR1174395}.}
Thus, Larson's theorem shows that even \textsf{MM} or \textsf{PFA} do not produce a Laver ideal over $\aleph_1$.
But in the main result of the current paper we are assuming that $2^\omega>\omega_2$.
This means that \textsf{MM} or \textsf{PFA} are not under the scope of our proof, and it leaves open the question of weakly Laver ideals over $\aleph_1$ under these forcing axioms.
Another aspect of this point is reflected in the attempt to deal with saturated ideals over larger accessible cardinals.
Dealing with such cardinals, we need to generalize Martin's axiom.

The discovery of Martin's axiom yielded a myriad of consistency results in many mathematical fields. The fact that one can add generic sets to all $ccc$ forcing notions simultaneously is very powerful. But the class of $ccc$ forcing notions is limited.

A natural question arose toward the end of the seventies: is it possible to generalize Martin's axiom to higher cardinals, with respect to other classes of forcing notions? The main objective was to phrase a parallel statement upon replacing $ccc$ by $\kappa$-cc forcing notions.

Three mathematicians worked out a positive answer almost at the same time. A version of generalized Martin's axiom has been phrased and proved by Shelah, Laver and Baumgartner.
Shelah's version appeared in \cite{MR0505492}, while Laver's version is an unpublished result.
Baumgartner's axiom comes from \cite{MR823775} and will be used in this paper.
All three variants are similar, and include additional requirements which do not appear in Martin's axiom. Later work showed that these additional features are obligatory in some sense.

Having a generalized version of Martin's axiom, one can try to prove statements which follow from Martin's axiom in the new setting. A comprehensive work in this line was carried out by Tall, \cite{MR1278025}.
We consider, in this paper, another issue.
We shall see that Larson's theorem generalizes to ideals over $\aleph_2$, upon replacing Martin's axiom by Baumgartner's axiom.

Actually, Larson's argument can be generalized in a broader context. Larson defined a cardinal characteristic called $\mathfrak{ap}$, and proved his theorem from the assumption that $\mathfrak{ap}=\omega_2=\mathfrak{c}$. Since this assumption follows from Martin's axiom and $2^\omega=\aleph_2$, the above theorem holds.

In recent years, an extensive literature is devoted to generalized cardinal characteristics. It follows that the definition of $\mathfrak{ap}$ applies equally well at every $\kappa=\cf(\kappa)$. Moreover, the pertinent assumption with respect to the generalized concept $\mathfrak{ap}_\kappa$ along with the non-existence of very dense ideals over $\kappa$ yield the corresponding conclusion about saturated ideals. The interesting question is whether a version of generalized Martin's axiom forces $\mathfrak{ap}_\kappa = 2^\kappa$.

By and large, not all the statements which follow from Martin's axiom can be generalized in this way. In particular, cardinal characteristics which assume the value of the continuum under Martin's axiom may be small under any version of generalized Martin's axiom. A salient example is the splitting number $\mathfrak{s}$ which equals $2^\omega$ under Martin's axiom. However, $\mathfrak{s}_{\aleph_1}=\aleph_0$ in \textsf{ZFC}. We shall see that the behavior of $\mathfrak{ap}_\kappa$ is similar to that of $\mathfrak{ap}$, and therefore the conclusion about saturated ideals can be generalized.

The argument of Larson can be implemented at higher cardinals, but the stronger theorem about weak Laverness proved in this paper seems more problematic.
Unlike Larson's proof, we do not know how to prove that Baumgartner's axiom excludes the existence of weakly Laver ideals over $\aleph_2$. This is pertinent, in particular, to the non-stationary ideal.

The paper contains three additional sections. The first one comprises of basic facts and required definitions.
In the second section we generalize Larson's proof and apply it to ideals over $\aleph_2$ under Baumgartner's axiom.
Finally, the third section contains the non-existence proof of weakly Laver ideals over $\aleph_1$ under Martin's axiom and $2^\omega>\aleph_2$.

We shall try to be consistent with notational conventions used in Larson's paper, \cite{MR2146222}. We use the Jerusalem forcing notation, so $p\leq q$ means that $p$ is weaker than $q$, compatible conditions are well met if they have a \emph{least upper bound}, a generic set is \emph{downward closed} and so on.
For general background concerning saturation properties of ideals and forcing with such ideals we refer to \cite{MR2768691} and \cite{MR2768692}.

\newpage

\section{Preliminaries}

All $\kappa$-ideals in this paper contain the ideal $J^{\rm bd}_\kappa$ of bounded subsets of $\kappa$. Likewise, all $\kappa$-ideals are $\kappa$-complete unless otherwise stated. If $\mathcal{I}$ is an ideal over $\kappa$ then $\mathcal{I}^+ = \mathcal{P}(\kappa)-\mathcal{I}$. The elements of $\mathcal{I}^+$ are called $\mathcal{I}$-positive sets. By the above conventions, if $\mathcal{I}$ is an ideal over $\kappa$ and $A\in\mathcal{I}^+$ then $|A|=\kappa$.
If $\mathcal{I}$ is a $\kappa$-ideal and $A,B\subseteq\kappa$ then $A\subseteq_{\mathcal{I}}B$ iff $A- B\in\mathcal{I}$.

\begin{definition}
\label{defsat} Saturation. \newline
An ideal $\mathcal{I}$ over $\kappa$ is $(\mu,\lambda,\theta)$-saturated iff for every collection of sets $\{A_\alpha:\alpha<\mu\}\subseteq\mathcal{I}^+$ there is a sub-collection $\{A_{\alpha_\beta}:\beta<\lambda\}$ such that for every $C\in[\lambda]^\theta$ it is true that $\bigcap\limits_{\beta\in C} A_{\alpha_\beta}\in\mathcal{I}^+$.
\end{definition}

Suppose that $\mathcal{I}$ is a $\kappa^+$-complete ideal over $\kappa^+$. The specific case of the above definition in which $\mu=\lambda=\kappa^{++}$ and $\theta=\kappa$ is the main subject of this paper. A $\kappa^+$-complete and $(\kappa^{++},\kappa^{++},\kappa)$-saturated ideal over $\kappa^+$ will be called Laver ideal. The existence of such an ideal over a successor of a regular cardinal has been established by Laver from a huge cardinal. It is unclear whether the existence of such an ideal is consistent over a successor of a singular cardinal. Such ideals yield many interesting results, some of which are connected with classical infinite combinatorics. There is, however, a stronger concept.

\begin{definition}
\label{defdense} Dense ideals. \newline
Let $\mathcal{J}$ be a $\kappa$-complete ideal over $\kappa$.
\begin{enumerate}
\item [$(\aleph)$] We say that $\mathcal{J}$ is $\nu$-dense iff there exists a collection $\{A_\alpha:\alpha<\nu\}\subseteq\mathcal{J}^+$ such that for every $B\in\mathcal{J}^+$ there is an ordinal $\alpha<\kappa$ for which $A_\alpha\subseteq_{\mathcal{J}}B$.
\item [$(\beth)$] The minimal $\nu$ for which $\mathcal{J}$ is $\nu$-dense is called the density of $\mathcal{J}$ and denoted by $d(\mathcal{J})$.
\item [$(\gimel)$] The weak density of $\mathcal{J}$, denoted by $wd(\mathcal{J})$, is the minimal $\nu$ such that there is a dense family $\mathcal{A} = \{A_\alpha:\alpha<\nu\}\subseteq \mathcal{J}^+$ with the property that every initial segment of $\mathcal{A}$ is not dense in $\mathcal{J}^+$.
\end{enumerate}
\end{definition}

It follows from the definition that $d(\mathcal{J})\leq wd(\mathcal{J})$, and this inequality will be used later.
The following simple observation shows that denseness is stronger than Laverness:

\begin{claim}
\label{clmdensity} A $\kappa^+$-dense ideal over $\kappa^+$ is Laver.
\end{claim}

\par\noindent\emph{Proof}. \newline
Let $\mathcal{J}$ be a $\kappa^+$-dense ideal over $\kappa^+$. Suppose that $\mathcal{B}=\{B_\delta:\delta<\kappa^{++}\}\subseteq\mathcal{J}^+$. Fix a dense collection $\{A_\alpha:\alpha<\kappa^+\}\subseteq\mathcal{J}^+$.
For every $\delta<\kappa^{++}$ choose an ordinal $\alpha(\delta)<\kappa^+$ such that $A_{\alpha(\delta)}\subseteq_{\mathcal{J}}B_\delta$. Since $\kappa^+ < \kappa^{++} = \cf(\kappa^{++})$, one can find a sub-collection $\{B_{\delta_\beta}:\beta<\kappa^{++}\}\subseteq\mathcal{B}$ and a fixed ordinal $\alpha<\kappa^+$ such that $\beta<\kappa^{++}\Rightarrow \alpha(\delta_\beta)=\alpha$.

We claim that $\{B_{\delta_\beta}:\beta<\kappa^{++}\}$ exemplifies the Laverness of $\mathcal{J}$ (with respect to $\mathcal{B}$). For this, assume that $\Gamma\subseteq\kappa^{++}, |\Gamma|\leq\kappa$. For every $\beta\in\Gamma$ we can see that $t_\beta=A_\alpha- B_{\delta_\beta}\in\mathcal{J}$. By the $\kappa^+$-completeness of $\mathcal{J}$ we see that $t\in\mathcal{J}$, so the set $A=A_\alpha- t\in\mathcal{J}^+$. Notice that $\bigcap\{B_{\delta_\beta}: \beta\in\Gamma\}\supseteq A$, and conclude that $\mathcal{J}$ is Laver.

\hfill \qedref{clmdensity}

We also mention the concept of $\lambda$-saturation. A $\kappa$-ideal $\mathcal{I}$ is $\lambda$-saturated iff for every collection $\{A_\alpha: \alpha<\lambda\}\subseteq\mathcal{I}^+$ one can find $\alpha<\beta<\lambda$ such that $A_\alpha\cap A_\beta\in\mathcal{I}^+$.
Thus, $\lambda$-saturation is the same as $(\lambda,2,2)$-saturation.
A uniform $\kappa^+$-complete $\kappa^{++}$-saturated ideal over $\kappa^+$ will be called Kunen ideal. Every Laver ideal is Kunen.

Let $\kappa$ be a regular cardinal. A family $\{e_\alpha:\alpha<\lambda\}\subseteq[\kappa]^\kappa$ is almost disjoint iff $\alpha<\beta<\lambda\Rightarrow |e_\alpha\cap e_\beta|<\kappa$. An almost disjoint family of size $2^\kappa$ always exists.
If $E=\{e_\alpha: \alpha<\lambda\}\subseteq[\kappa]^\kappa$ is almost disjoint and $A\subseteq\lambda$ then a set $x\subseteq\kappa$ \emph{depicts} the pair $(A,E)$ iff $|e_\alpha\cap x|=\kappa\Leftrightarrow\alpha\in A$. If there is no $x\subseteq\kappa$ which depicts $(A,E)$ then we shall say that the pair $(A,E)$ is \emph{indepictable}.

\begin{definition}
\label{defapkappa} $\mathfrak{ap}_\kappa$. \newline
Let $\kappa$ be a regular cardinal. \newline
The characteristic $\mathfrak{ap}_\kappa$ is the minimal cardinal $\lambda$ for which there is an almost disjoint family $E=\{e_\alpha: \alpha<\lambda\}\subseteq[\kappa]^\kappa$ and a subset $A\subseteq\lambda$ such that $(A,E)$ is indepictable.
\end{definition}

An immediate question is whether $\mathfrak{ap}_\kappa$ is well defined.
The following claim gives a positive answer, using a simple property of this characteristic.

\begin{claim}
\label{clmapwell} $\mathfrak{ap}_\kappa$ is well-defined whenever $\kappa = \cf(\kappa)$. Likewise, if $\kappa\leq\theta<\mathfrak{ap}_\kappa$ then $2^\theta=2^\kappa$.
\end{claim}

\par\noindent\emph{Proof}. \newline
Assume that $E=\{e_\alpha: \alpha<\lambda\}\subseteq[\kappa]^\kappa$ is almost disjoint. Let $A,B\subseteq\lambda$ so that $A\neq B$. We observe that if $x(A)$ depicts $(A,E)$ and $x(B)$ depicts $(B,E)$ then $x(A)\neq x(B)$. Indeed, $A\neq B$ so without loss of generality there exists an ordinal $\alpha\in A- B$. It follows that $|e_\alpha\cap x(A)|=\kappa$ while $|e_\alpha\cap x(B)|<\kappa$, so in particular $x(A)\neq x(B)$.

Now let $\lambda=2^\kappa$ and fix any almost disjoint family $E=\{e_\alpha: \alpha<\lambda\}\subseteq[\kappa]^\kappa$. For every $A\subseteq\lambda$ choose $x(A)\subseteq\kappa$ such that $x(A)$ depicts $(A,E)$, if possible. By this, a one-to-one partial mapping is created. However, $|\mathcal{P}(\lambda)|>\lambda = |\mathcal{P}(\kappa)|$, so necessarily there is a subset $A\subseteq\lambda$ such that $(A,E)$ is indepictable. It follows that $\mathfrak{ap}_\kappa$ is well-defined (and of course, $\mathfrak{ap}_\kappa\leq 2^\kappa$). The above argument shows that $2^\theta=2^\kappa$ whenever $\theta\in[\kappa,\mathfrak{ap}_\kappa)$, so we are done.

\hfill \qedref{clmapwell}

Let us try to explain the idea behind the definition of $\mathfrak{ap}_\kappa$.
If $e\in[\kappa]^\kappa$ and $E=\{e\}$ then $E$ is almost disjoint. But this trivial example gives no information about the structure of $[\kappa]^\kappa$. In order to be more informative we need a \emph{rich} almost disjoint family.

How are we to gauge richness? One way is to ask for a maximal almost disjoint family. The minimal cardinality of such a family is the traditional almost disjointness number $\mathfrak{a}_\kappa$.
The characteristic $\mathfrak{ap}_\kappa$ gives another interpretation for richness.
Notice that if $E$ is maximal almost disjoint in $[\kappa]^\kappa$ then $(\varnothing,E)$ is indepictable, and hence $\mathfrak{ap}_\kappa\leq\mathfrak{a}_\kappa$.

We shall need a description of Baumgartner's axiom. Let $\mathbb{P}$ be a forcing notion and $L\subseteq\mathbb{P}$. We say that $L$ is linked iff the conditions in $L$ are pairwise compatible. The forcing notion $\mathbb{P}$ is $\aleph_1$-linked iff $\mathbb{P}$ is expressible as $\bigcup_{\gamma\in\omega_1}L_\gamma$ when each $L_\gamma$ is linked. A forcing notion $\mathbb{P}$ is well-met iff any pair of compatible elements has a least upper bound.

\begin{axiom}
\label{axba} Baumgartner's axiom. \newline
For every countably closed $\aleph_1$-linked well met $\mathbb{P}$ and every collection $\{D_\alpha:\alpha<\theta\}$ of dense subsets of $\mathbb{P}$ such that $\theta<2^{\aleph_1}$, there exists a generic set $G\subseteq\mathbb{P}$ such that $G\cap D_\alpha\neq\varnothing$ for every $\alpha<\theta$.
\end{axiom}

Baumgartner proved in \cite{MR823775} the following:

\begin{theorem}
\label{thmba} Assume that $\kappa=\cf(\kappa)>\aleph_1$ and $2^{<\kappa}=\kappa$. Then one can force Baumgartner's axiom with $2^{\aleph_1}=\kappa$. If $2^{\aleph_0}<\kappa$ in the ground model then one can force $2^{\aleph_0}=\aleph_1, 2^{\aleph_1}=\kappa$ and Baumgartner's axiom.
\end{theorem}

\hfill \qedref{thmba}

We conclude this section with the definition of polarized partition relations. This will be the main tool for handling weakly Laver ideals in the last section.

\begin{definition}
\label{defpolarized} Polarized partition relations. \newline
The relation $\binom{\lambda}{\kappa} \rightarrow \binom{\mu}{\nu}^{1,1}_\theta$ means that for every coloring $c:\lambda\times\kappa\rightarrow\theta$ one can find a pair of sets $A\subseteq\lambda, B\subseteq\kappa$ so that $|A|=\mu, |B|=\nu$ and $c\upharpoonright(A\times B)$ is constant.
\end{definition}

\newpage

\section{Laverness}

We open this section with a description of $\mathfrak{ap}_{\aleph_1}$ under Baumgartner's axiom. Recall that under Martin's axiom the value of $\mathfrak{ap}$ is the continuum, but this is not surprising since almost every cardinal characteristic known to man equals the continuum under Martin's axiom.\footnote{An outlier is the Galvin number, defined in \cite{MR3787522}. Consistently, this cardinal characteristic is smaller than the continuum under Martin's axiom. } This is not the case with respect to Baumgartner's axiom. However, it holds for $\mathfrak{ap}_{\aleph_1}$, as indicated in the following theorem.

\begin{theorem}
\label{thmbaumap} Assume Baumgartner's axiom, $2^{\aleph_0}=\aleph_1$ and $2^{\aleph_1}=\kappa$. \newline
Then $\mathfrak{ap}_{\aleph_1} = \kappa$.
\end{theorem}

\par\noindent\emph{Proof}. \newline
Suppose that $\theta<\kappa$, and assume toward contradiction that $\mathfrak{ap}_{\aleph_1} = \theta$. Choose an almost disjoint family $E=\{e_\alpha:\alpha<\theta\}\subseteq [\aleph_1]^{\aleph_1}$ and $A\subseteq\theta$ such that $(A,E)$ is indepictable. We define a forcing notion that adds a set $x\subseteq\omega_1$ which depicts $(A,E)$, thus arriving at a contradiction.

A condition $p\in\mathbb{P}$ is a pair $(t,B)=(t_p,B_p)$, where $t\subseteq\omega_1, B\subseteq\theta- A$ and $|t|=|B|=\aleph_0$. If $p,q\in\mathbb{P}$ then $p\leq_{\mathbb{P}}q$ iff $t_p\subseteq t_q, B_p\subseteq B_q$ and if $\gamma\in t_q- t_p$ then for every $\alpha\in B_p$ one has $\gamma\notin e_\alpha$. The left-hand coordinate $t_p$ approximates a set $x\subseteq\omega_1$. This set will depict the pair $(A,E)$ in the generic extension. The right-hand coordinate $B_p$ gives some information about elements $e_\alpha\in{E}$, with respect to every $\alpha\in B_p$. This information makes sure that $A$ is depictable in the generic extension.

In order to claim for the existence of a generic object for $\mathbb{P}$ we must show that it satisfies the requirements of Baumgartner's axiom. First we notice that $\mathbb{P}$ is countably closed. Indeed, assume that $\langle p_n: n\in\omega\rangle$ is an increasing sequence of conditions. Define $t=\bigcup_{n\in\omega}t_{p_n}, B=\bigcup_{n\in\omega}B_{p_n}$ and $q=(t,B)$. It follows from the definition of the forcing order that $\forall n\in\omega, p_n\leq q\in\mathbb{P}$, so $\mathbb{P}$ is countably closed.

Assume now that $p=(t,B_0)$ and $q=(t,B_1)$. If $r=(t,B_0\cup B_1)$ then $p,q\leq r\in\mathbb{P}$, so any collection of conditions with the same left coordinate is linked. Recall that $2^{\aleph_0}=\aleph_1$, so $|\{t: t\subseteq\omega_1, |t|=\aleph_0\}|=\aleph_1$ as well. For every $t\in[\omega_1]^{\aleph_0}$ let $L_t = \{p\in\mathbb{P}: t_p=t\}$. Each $L_t$ is linked, and since $\mathbb{P} = \bigcup\{L_t: t\in[\omega_1]^{\aleph_0}\}$ we conclude that $\mathbb{P}$ is $\aleph_1$-linked.

Finally, suppose that $p\parallel q$. Let $t=t_p\cup t_q, B=B_p\cup B_q$ and $r=(t,B)$. Since $p$ and $q$ are compatible, there is a condition $s\in\mathbb{P}$ so that $s\geq p,q$. By the definition of the forcing order, $t_s\supseteq t_p,t_q$ and $B_s\supseteq B_p,B_q$. This shows that $r\in\mathbb{P}$, so in particular it is an upper bound of $p,q$. Moreover, it shows that $r$ is a least upper bound, so $\mathbb{P}$ is well-met.

For every $\alpha\in\theta- A$ let $D_\alpha = \{p\in\mathbb{P}: \alpha\in B_p\}$. Fix an ordinal $\alpha\in\theta- A$, and assume that $q\notin D_\alpha$. Define $r=(t_q,B_q\cup\{\alpha\})$. Clearly $q\leq r\in D_\alpha$ so each $D_\alpha$ is a dense subset of $\mathbb{P}$.
For every $\alpha\in A$ and each $\beta\in\omega_1$ let $E_{\alpha\beta} = \{p\in\mathbb{P}: \exists\gamma\geq\beta, \gamma\in e_\alpha\cap t_p\}$.
Fix $\alpha\in{A},\beta\in\omega_1$.
Assume that $q\notin E_{\alpha\beta}$.
Notice that $\alpha\notin B_q$ as $\alpha\in A, B_q\subseteq\theta- A$, and hence for every $\alpha_i\in B_q$ we have $|e_\alpha\cap e_{\alpha_i}|<\aleph_1$. Consequently, one can choose an ordinal $\gamma\geq\beta$ such that $\gamma\in e_\alpha - \bigcup_{\alpha_i\in B_q} e_{\alpha_i}$. Define $r=(t_q\cup\{\gamma\},B_q)$. It follows that $q\leq r\in E_{\alpha\beta}$ and hence each $E_{\alpha\beta}$ is dense.

Let $\mathcal{D}=\{D_\alpha\mid\alpha\in\theta-A\}\cup\{E_{\alpha\beta}\mid\alpha\in{A},\beta\in\omega_1\}$.
Observe that $|\mathcal{D}|\leq\theta<\kappa$, and choose a $V$-generic set $G\subseteq\mathbb{P}$ which meets every element of $\mathcal{D}$.
Let $x=x_G = \bigcup\{t: \exists B, (t,B)\in G\}$.
Let $\name{x}$ be a canonical name for $x$.
It follows that $|x\cap e_\alpha|<\aleph_1$ for every $\alpha\in\theta-A$ since $G\cap D_\alpha\neq\varnothing$, and $|x\cap e_\alpha|=\aleph_1$ for every $\alpha\in{A}$ since $G\cap E_{\alpha\beta}\neq\varnothing$ for every $\beta\in\omega_1$.
Thus the pair $(A,E)$ is depictable, a contradiction.

\hfill \qedref{thmbaumap}

Our second theorem deals with dense ideals over $\aleph_2$.
In the theorem below we focus on dense ideals with an additional closure property.
Call $\mathcal{I}$ strongly $\omega_1$-closed if every $\subseteq_{\mathcal{I}}$-decreasing sequence $(A_i:i\in\omega)$ of elements of $\mathcal{I}^+$ admits a lower bound, namely $A\in\mathcal{I}^+$ so that $A\subseteq_{\mathcal{I}}A_i$ for every $i\in\omega$.
We shall prove that under Baumgartner's axiom (and appropriate cardinal arithmetic) there are no strongly $\omega_1$-closed $\aleph_2$-dense $\aleph_2$-complete ideals.
This is also the parallel situation under Martin's axiom with respect to $\aleph_1$, see e.g. \cite[Theorem 7.6]{MR530430}.
However, we shall take a different path than the one in Taylor's proof.
The gist of the proof below is exactly as in the proof of Theorem \ref{thmbaumap} (actually, it seems that the statement will follow merely from $\mathfrak{ap}_{\aleph_1}>\aleph_2$, which holds under Baumgartner's axiom as proved above). It is a generalization of the parallel argument in \cite{MR924672} under Martin's axiom.

\begin{theorem}
\label{thmdense} Assume $\mathfrak{ap}_{\aleph_1}>\aleph_2$.
Then there is no strongly $\omega_1$-closed $\aleph_2$-dense $\aleph_2$-complete ideal over $\aleph_2$.
Hence if Baumgartner's axiom holds and $2^{\aleph_0}=\aleph_1, 2^{\aleph_1}>\aleph_2$, then there are no such ideals over $\aleph_2$.
\end{theorem}

\par\noindent\emph{Proof}. \newline
Suppose that $\mathcal{J}$ is an $\omega_1$-strongly closed $\aleph_2$-dense ideal over $\aleph_2$. Let $\mathbb{Q}$ be the forcing notion whose conditions are the elments of $\mathcal{J}^+$ ordered by reversed inclusion. This forcing notion adds a Cohen $\omega_1$-real (i.e., an element of ${}^{\omega_1}2$), which is not minimal over the ground model $V$.

Let us explain why a Cohen subset is added.
We describe, firstly, the parallel situation at $\aleph_1$.
Let $\mathcal{I}$ be an $\aleph_1$-dense ideal over $\aleph_1$, and let $\mathbb{Q}$ be the associated forcing (to wit, the conditions of $\mathbb{Q}$ are the elements of $\mathcal{I}^+$ mod $\mathcal{I}$).
Then $\mathcal{P}(\omega_1)/\mathcal{I}\cong{\rm Col}(\omega,\omega_1)$ and hence forcing with $\mathbb{Q}$ adds a Cohen real (see \cite{MR2768692}, or \cite[Corollary 19]{MR924672}).

Essentially, the same is true with respect to $\aleph_2$-dense ideals over $\aleph_2$, but here we need the additional assumption that the corresponding forcing notion $\mathbb{Q}$ is $\aleph_1$-closed (that is, $\mathcal{I}$ is $\omega_1$-strongly closed).
In this setting, $\mathcal{P}(\omega_2)/\mathcal{I}\cong{\rm Col}(\omega_1,\omega_2)$ and one can build a generic for ${\rm Add}(\omega_1,1)$ by induction on $\omega_1$.
However, in this case one has to overcome countable sequences along the process, and the assumption that $\mathbb{Q}$ is $\aleph_1$-closed comes in handy.
These phenomena are elaborated in \cite{MR2768691}.

Assume now that Baumgartner's axiom holds, $2^{\aleph_0}=\aleph_1$ and $2^{\aleph_1}>\aleph_2$. Let $\mathcal{I}$ be an $\aleph_2$-dense $\aleph_2$-complete ideal over $\aleph_2$.
It follows that $\mathcal{I}$ is precipitous (actually, $\omega_3$-saturation suffices).
Let $G\subseteq \mathcal{P}(\omega_2)/\mathcal{I}$ be generic. Assume that $r$ is an $\omega_1$-real, $r\in V[G]$ but $r\notin V$. We shall see that $r$ is minimal over $V$, i.e. $V[r]=V[G]$.
In other words, every subset of $\omega_1$ in $V[G]-V$ is minimal over $V$.
In the light of the above paragraph, this concludes the proof of the theorem.

Let $M$ be the transitive collapse of $V^{\omega_2}/G$, and let $\jmath:V\rightarrow M$ be the generic embedding.
Notice that $({}^{\omega_1}2)^M = ({}^{\omega_1}2)^{V[G]}$, so $r\in M$ and hence represented by an element in $V$.
As a reference to the above fact we suggest \cite[Section 1.6]{MR2069032}.
The theorems there are phrased at the level of $\aleph_1$, but similar arguments work at $\aleph_2$, upon noticing that an $\aleph_2$-dense ideal over $\aleph_2$ is saturated.
Fix a function $f:\omega_2\rightarrow{}^{\omega_1}2$ so that $f\in V$ and $r=[f]_M$. Without loss of generality, $f$ is one-to-one.
Enumerate the elements of ${}^{<\omega_1}2$ by $\{\tau_\eta:\eta\in\omega_1\}$ (here we use the assumption that $2^\omega=\omega_1$). For every $s\in {}^{\omega_1}2$ let $e_s = \{\eta\in\omega_1: \exists\zeta\in\omega_1, s\upharpoonright\zeta=\tau_\eta\}$.
Notice that $\{e_s:s\in{}^{\omega_1}2\}$ forms an almost disjoint family.

We claim that if $A\subseteq\omega_2$ then for some $x\subseteq\omega_1$ we have $|x\cap e_{f(\alpha)}|\leq\aleph_0$ iff $\alpha\in A$. In the phraseology of cardinal characteristics it means exactly that every $A\subseteq\omega_2$ is depictable with respect to the almost disjoint collection of the $e_{f(\alpha)}$-s. This statement is actually a consequence of the assumption that $\mathfrak{ap}_{\aleph_1}>\aleph_2$, and hence follows from Baumgartner's axiom as proved in Theorem \ref{thmbaumap}.

But now, if $A\subseteq\omega_2, A\in{V}$ and $x\subseteq\omega_1, x\in{V}$ depicts $A$ in the sense that $\alpha\in A$ iff $|x\cap e_{f(\alpha)}|\leq\aleph_0$, then $A\in G$ iff $e_r\cap x$ is countable. Since $\mathcal{P}(\omega_2)/\mathcal{I}$ is $\sigma$-complete, no countable subset of the ground model is added in the extension. Consequently, $G$ can be recovered from $r$ by collecting all the countable sets of the form $e_r\cap x$. It follows that $V[r]=V[G]$, as required.

\hfill \qedref{thmdense}

The rest of this section generalizes Larson's proof. As Larson indicates for the case of $\aleph_1$, the proof is purely combinatorial and follows from the fact that $\mathfrak{ap}_\kappa$ assumes a large value (and the additional fact that there are no $\kappa$-dense ideals over $\kappa$).

\begin{definition}
\label{defenx} Let $\{e_\gamma:\gamma\in\omega_2\}\subseteq [\omega_1]^{\omega_1}$ be an almost disjoint family, $x\subseteq\omega_1$ and $\eta\in\omega_1$. \newline
We let $E^\eta_x = \{\gamma\in\omega_2:{\rm otp}(e_\gamma\cap x)\geq\eta\}$ and $F^\eta_x = \omega_2- E^\eta_x$.
\end{definition}

The following lemma provides a pair of large disjoint sets $(E^\eta_x, F^\eta_x)$ under some circumstances:
We indicate that an almost disjoint family as in Definition \ref{defenx} is fixed, and the sets $E^\eta_x,F^\eta_x$ refer to this family.

\begin{lemma}
\label{lemfnx} Assume that:
\begin{enumerate}
\item [$(a)$] $\mathfrak{ap}_{\aleph_1}>\aleph_2$.
\item [$(b)$] $\mathcal{I}$ is an $\aleph_2$-complete ideal over $\aleph_2$.
\item [$(c)$] $\mathcal{P}(\omega_2)/\mathcal{I}$ is not $\theta$-dense.
\item [$(d)$] $\{A_\alpha:\alpha<\theta\}\subseteq\mathcal{I}^+$.
\end{enumerate}
Then one can find $x\subseteq\omega_1$ and $\eta\in\omega_1$ for which $F^\eta_x\in\mathcal{I}^+$ and for every $\alpha<\theta$ there is some $\zeta\in\omega_1$ such that $E^\eta_{x\cap\zeta}\cap A_\alpha\in\mathcal{I}^+$.
\end{lemma}

\par\noindent\emph{Proof}. \newline
Fix an almost disjoint family $E=\{e_\gamma:\gamma\in\omega_2\}$ in $[\omega_2]^{\omega_2}$. As a first step we refine the collection $\{A_\alpha: \alpha<\theta\}$ of $\mathcal{I}$-positive sets in the following way. We shall construct a collection $\{B_\alpha:\alpha<\theta\}\subseteq\mathcal{I}^+$ and a set $D\in\mathcal{I}^+$ such that $\forall\alpha<\theta, B_\alpha\subseteq A_\alpha$ and $B_\alpha\cap D=\varnothing$.

By assumption $(c)$, the collection $\{A_\alpha:\alpha<\theta\}$ is not dense in $\mathcal{P}(\omega_2)/\mathcal{I}$. Hence there is an element $D\in\mathcal{I}^+$ such that $\forall\alpha\in\theta, \neg(A_\alpha \subseteq_{\mathcal{I}}D)$, i.e. $B_\alpha = A_\alpha- D\in\mathcal{I}^+$. Notice that $\{B_\alpha: \alpha<\theta\}$ is as required.

By the construction, $\bigcup_{\alpha<\theta}B_\alpha\cap D=\varnothing$. Now $\mathfrak{ap}_{\aleph_1}>\aleph_2$ and hence every subset of $\omega_2$ is depictable with $E$. In particular, the pair $(\omega_2- D,E)$ is depictable, so one can find $x\subseteq\omega_1$ such that:
\begin{itemize}
\item $|e_\gamma\cap x|<\aleph_1$ for every $\gamma\in D$.
\item $|e_\gamma\cap x|=\aleph_1$ for every $\gamma\notin D$.
\end{itemize}
It follows that $D\subseteq\bigcup\{F^\eta_x:\eta\in\omega_1\}$. For proving this fact, fix an ordinal $\gamma\in D$. Since $|e_\gamma\cap x|<\aleph_1$ we can choose an ordinal $\eta\in\omega_1$ such that ${\rm otp}(e_\gamma\cap x)<\eta$. By definition $\gamma\in F^\eta_x$, so $D\subseteq\bigcup\{F^\eta_x:\eta\in\omega_1\}$ indeed.

Since $\mathcal{I}$ is $\aleph_2$-complete and $D\in\mathcal{I}^+$ we infer that $F^\eta_x\notin\mathcal{I}$ for some $\eta\in\omega_1$. So we fix such an ordinal $\eta$, and the corresponding set $F^\eta_x$ will satisfy the first statement of the lemma. We must show that the same $x$ and $\eta$ satisfy the second statement as well. So assume that $\alpha<\theta$. We claim that $B_\alpha\subseteq\bigcup\{E^\eta_{x\cap\zeta}:\zeta\in\omega_1\}$. This follows since $B_\alpha\cap D = \varnothing$, so if $\gamma\in B_\alpha$ then $|e_\gamma\cap x|=\aleph_1$. Hence for some $\zeta\in\omega_1$ we must have ${\rm otp}(e_\gamma\cap(x\cap\zeta))\geq\eta$ and then $\gamma\in E^\eta_{x\cap\zeta}$. But this means that $E^\eta_{x\cap\zeta}\cap B_\alpha\in\mathcal{I}^+$ for some $\zeta\in\omega_1$, since $B_\alpha\in\mathcal{I}^+$ and $\mathcal{I}$ is $\aleph_2$-complete. Consequently, $E^\eta_{x\cap\zeta}\cap A_\alpha\in\mathcal{I}^+$ as $B_\alpha\subseteq A_\alpha$, so we are done.

\hfill \qedref{lemfnx}

The role of the lemma can be described as follows. Any application of the lemma gives ``one-step" of non-saturation. The set $x\subseteq\omega_1$ ensured by the lemma gives rise to a pair of \emph{big disjoint sets}. We shall prove non-saturation by an inductive process in which at every step we employ the lemma to produce more and more such pairs.

\begin{theorem}
\label{thmmt} Assume that:
\begin{enumerate}
\item [$(a)$] $\mathfrak{ap}_{\aleph_1}>\aleph_2$.
\item [$(b)$] $\mathcal{I}$ is an $\aleph_2$-complete ideal over $\aleph_2$.
\item [$(c)$] $wd(\mathcal{I})=\theta$.
\item [$(d)$] $2^\omega=\omega_1$.
\end{enumerate}
Then one can find a sequence $\langle D_\beta:\beta<\theta\rangle$, each $D_\beta$ is an element of $\mathcal{I}^+$, such that whenever $c\subseteq\theta$ is unbounded in $\theta$ there is a countable $y\subseteq c$ for which $\bigcap\{D_\beta:\beta\in y\}\in\mathcal{I}$.
\end{theorem}

\par\noindent\emph{Proof}. \newline
Let $\{e_\gamma:\gamma\in\omega_2\}\subseteq[\omega_1]^{\omega_1}$ be almost disjoint. The sets $E^\eta_x, F^\eta_x$ from Definition \ref{defenx} are defined with respect to this family. Fix an $\mathcal{I}^+$-dense collection $\{A_\alpha:\alpha<\theta\}\subseteq \mathcal{I}^+$ such that the initial segment $\{A_\alpha:\alpha<\beta\}$ is not $\mathcal{I}^+$-dense for every $\beta<\theta$.

By Lemma \ref{lemfnx} we choose, for each $\beta<\theta$, a set $x_\beta\subseteq\omega_1$ and an ordinal $\eta_\beta\in\omega_1$ so that $F^{\eta_\beta}_{x_\beta}\in\mathcal{I}^+$ and for every $\alpha<\beta$ there is some $\zeta_\alpha=\zeta\in\omega_1$ such that $E^{\eta_\beta}_{x_\beta\cap \zeta}\cap A_\alpha\in\mathcal{I}^+$. We denote $F^{\eta_\beta}_{x_\beta}$ by $D_\beta$ for every $\beta<\theta$. We claim that $\langle D_\beta: \beta\in\theta\rangle$ satisfies the theorem. So assume that $c\subseteq\theta$ is unbounded, and define:
$$
Z = \{(\eta,\sigma): \eta\in\omega_1, \sigma\in[\omega_1]^{<\omega_1}, \exists\beta\in c (E^\eta_\sigma\cap D_\beta\in\mathcal{I})\}.
$$
Notice that $|Z|=\aleph_1$ by assumption $(d)$ of the theorem. The crucial point of the proof is the fact that for every $\alpha<\theta$ one can find $\eta\in\omega_1, \sigma\in[\omega_1]^{<\omega_1}$ and $\beta\in c$ such that $E^\eta_\sigma\cap D_\beta\in\mathcal{I}$ while $E^\eta_\sigma\cap A_\alpha\in \mathcal{I}^+$.

For proving the crucial point, fix an ordinal $\alpha<\theta$. Since $c$ is unbounded in $\theta$ we can pick up some $\beta\in c$ so that $\beta>\alpha$. Lemma \ref{lemfnx} gives now $\eta_\beta = \eta\in\omega_1, x_\beta = x\subseteq\omega_1$ and $\zeta\in\omega_1$ such that $E^\eta_{x\cap\zeta}\cap A_\alpha\in\mathcal{I}^+$. Let $\sigma = x\cap\zeta$. We claim that the triple $(\eta,\sigma,\beta)$ exemplifies the crucial point with respect to $\alpha$.
Indeed, $E^\eta_\sigma\cap A_\alpha = E^\eta_{x\cap\zeta}\cap A_\alpha\in \mathcal{I}^+$ by the choice of $(x,\eta,\zeta)$. Likewise, $D_\beta = F^{\eta_\beta}_{x_\beta} = F^\eta_x$. Now $E^\eta_\sigma\cap F^\eta_x = E^\eta_\sigma\cap D_\beta \in\mathcal{I}$.

Having the crucial point at hand, we can finish the proof. For every pair $(\eta,\sigma)\in Z$ we choose an ordinal $\beta_{\eta\sigma}\in c$ so that $E^\eta_\sigma\cap D_{\beta_{\eta\sigma}}\in\mathcal{I}$ (such an ordinal exists by the definition of $Z$). Let $y = \{\beta_{\eta\sigma}:(\eta,\sigma)\in Z\}$. Notice that $y\subseteq c$ and $|y|\leq\aleph_1$. Define $y' = \{D_{\beta_{\eta\sigma}}:\beta_{\eta\sigma}\in y\}$. It remains to show that $\bigcap y'\in\mathcal{I}$.

For this end notice that if $(\eta,\sigma)\in Z$ then $\bigcap y'\cap E^\eta_\sigma\subseteq D_{\beta_{\eta\sigma}}\cap E^\eta_\sigma$. Assume toward contradiction that $\bigcap y'\in\mathcal{I}^+$. By density, we can find an ordinal $\alpha<\theta$ such that $A_\alpha\subseteq_{\mathcal{I}}\bigcap y'$. Now we can find $(\eta,\sigma)\in Z$ for which $E^\eta_\sigma\cap A_\alpha\in \mathcal{I}^+$ and conclude that $\bigcap y'\cap E^\eta_\sigma\in \mathcal{I}^+$, a contradiction.

To sum up, the collection $\langle D_\beta:\beta<\theta\rangle$ satisfies the property that for each unbounded $c\subseteq\theta$ there is a subset $y\subseteq c, |y|\leq\aleph_1$ such that $\bigcap\{D_\beta:\beta\in y\}\in\mathcal{I}$, so the proof is accomplished.

\hfill \qedref{thmmt}

Combining all the above facts, we have:

\begin{corollary}
\label{cormcor} Assume Baumgartner's axiom and $2^{\aleph_1}\geq\aleph_3$. \newline
Then there is no strongly $\omega_1$-closed Laver ideal over $\aleph_2$.
Moreover, there is no strongly $\omega_1$-closed $(\aleph_3,\aleph_3,\aleph_0)$-saturated ideal over $\aleph_2$.
\end{corollary}

\par\noindent\emph{Proof}. \newline
Let $\mathcal{I}$ be an $\aleph_2$-complete ideal over $\aleph_2$. From Theorem \ref{thmbaumap} we know that $\mathfrak{ap}_{\aleph_1}=2^{\aleph_1}$. From Theorem \ref{thmdense} we know that $d(\mathcal{I})\geq\aleph_3$. Hence Theorem \ref{thmmt} applies, and its conclusion shows that $\mathcal{I}$ is not Laver.

\hfill \qedref{cormcor}

We conclude this section with several comments. First, the proofs under Baumgartner's axiom can be implemented with Shelah's version of the generalized Martin's axiom from \cite{MR0505492}. Second, the case of $\aleph_2$ is only one example, and the theorem is general. For every $\kappa$ such that Baumgartner's axiom (or Shelah's generalized Martin's axiom) applies, one can prove the non-existence of Laver ideals at $\kappa^+$. The main point is that the natural forcing poset which increases $\mathfrak{ap}_\kappa$ behaves nicely.

Finally, it would be interesting to know whether the assumption that $\mathcal{I}$ is strongly $\omega_1$-closed can be removed from the above statements.
Observe that in order to carry out the argument within the proof of Theorem \ref{thmdense} one has to show that forcing with $\mathcal{I}^+$ results in a non-minimal generic extension.
By introducing a Cohen subset one obtains this goal, but it is plausible that a non-minimal generic extension (in this setting) is created even if forcing with $\mathcal{I}^+$ does not add a Cohen subset.

\begin{question}
  \label{qminimalextension} Let $\mathbb{Q}$ be a forcing notion which collapses $\aleph_2$ to $\aleph_1$, and let $G\subseteq\mathbb{Q}$ be $V$-generic. Is it possible that $V[A]=V[G]$ whenever $A\subseteq\omega_1$ and $A\notin{V}$?
\end{question}

\newpage

\section{Echeloned saturation}

In this section we address the question of Larson about the existence of echeloned ideals over $\aleph_1$ under forcing axioms.
The main result reads as follows.
Suppose that Chang's conjecture fails, and assume that Martin's axiom holds and $2^{\aleph_0}>\aleph_2$.
Then there is no weakly Laver ideal (i.e., an echeloned ideal) over $\aleph_1$.
Let us commence with the formal definition of weak Laverness.

\begin{definition}
\label{defweak} Weakly Laver ideals. \newline
An ideal $\mathcal{I}$ over a successor cardinal $\kappa^+$ is weakly Laver if $\mathcal{I}$ is uniform, $\kappa^+$-complete and $(\kappa^{++},\kappa^+,\kappa)$-saturated.
\end{definition}

Larson, \cite{MR2146222}, asked whether ${\rm NS}_{\omega_1}$ is $(\aleph_2,\aleph_1,\aleph_0)$-saturated under some forcing axiom. Being uniform and $\sigma$-complete, this amounts to the question whether ${\rm NS}_{\omega_1}$ is weakly Laver. We shall give a negative answer, not only with respect to ${\rm NS}_{\omega_1}$ but for every uniform $\sigma$-complete ideal over $\aleph_1$, provided that Chang's conjecture fails.

Our strategy is as follows.
On the one hand, if Chang's conjecture fails then ${\rm MA}_{\aleph_2}$ implies a negative combinatorial relation.
This is a result of Todor\v{c}evi\'c.
On the other hand, the existence of weakly Laver ideals implies a positive relation with the same parameters.
Thus, under the failure of Chang's conjecture, ${\rm MA}_{\aleph_2}$ shows that there is no $(\aleph_2,\aleph_1,\aleph_0)$-saturated ideal over $\aleph_1$.
We indicate that the failure of Chang's conjecture is, in some sense, necessary if one wishes to pursue this strategy, as we shall explain below.

Recall that Chang's conjecture is the statement that for every structure $\mathfrak{A}=(A,R,\ldots)$ in a first order countable language with $|A|=\aleph_2$ and $|R|=\aleph_1$ there exists $\mathfrak{B}\prec\mathfrak{A}$ such that $|B|=\aleph_1$ and $|R^\mathfrak{B}|=\aleph_0$.
The following appears in \cite[Theorem 3]{MR1127033}.

\begin{theorem}
  \label{thmtodorcevic} Assume ${\rm MA}_{\aleph_2}$.
  Then Chang's conjecture is equivalent to the positive relation $\binom{\omega_2}{\omega_2}\rightarrow\binom{\omega}{\omega}^{1,1}_\omega$.
\end{theorem}

\hfill \qedref{thmtodorcevic}

Roughly, the idea of Todor\v{c}evi\'c is based on a characterization of the failure of Chang's conjecture by the existence of a function $a:[\omega_2]^2\rightarrow\omega_1$ with some prescribed properties \cite[Lemma 4]{MR1127033}.
Todor\v{c}evi\'c defines a forcing notion which implies $\binom{\omega_2}{\omega_2}\nrightarrow\binom{\omega}{\omega}^{1,1}_\omega$, and proves that if the above function $a$ exists then this forcing notion is $ccc$.
Moreover, one needs just $\aleph_2$-many dense subsets in order to carry out the argument.

\begin{corollary}
  \label{cortodorcevic} Assume the failure of Chang's conjecture.
  If Martin's axiom holds and $2^{\aleph_0}>\aleph_2$ then $\binom{\omega_2}{\omega_1}\nrightarrow\binom{\omega}{\omega}^{1,1}_\omega$.
\end{corollary}

\par\noindent\emph{Proof}. \newline
From the above theorem we see that $\binom{\omega_2}{\omega_2}\nrightarrow\binom{\omega}{\omega}^{1,1}_\omega$.
Now use monotonicity.

\hfill \qedref{cortodorcevic}

Let us mention the fact that in some sense the failure of Chang's conjecture is necessary here.
Indeed, Todor\v{c}evi\'c proved in \cite{MR1127033} that Chang's conjecture implies $\binom{\omega_2}{\omega_2}\rightarrow\binom{\omega}{\omega}^{1,1}_\omega$.
Inasmuch as Chang's conjecture is indestructible under $ccc$ forcing notions and Martin's axiom can be forced by a $ccc$ forcing notion, one concludes that the negation of Chang's conjecture is essential in our argument.

Concerning the above theorem we indicate that the number of colors, $\aleph_0$, is also essential for the negative relation that has been proved. If one replaces it by any finite number of colors then strong positive relations hold under Martin's axiom.

We shall prove now the opposite assertion, using weakly Laver ideals over $\aleph_1$. The proof is modelled after \cite{MR3610266}, in which a much stronger theorem is proved under the existence assumption of Laver ideals. The main point of departure is the fact that no special closure properties are assumed with respect to the elementary sub-models that we shall use here. It gives a weaker combinatorial result, but it enables us to prove the theorem without assuming that the continuum is small.

\begin{theorem}
\label{thmweaklaver} Assume that there exists a weakly Laver ideal over $\aleph_1$. \newline
Then $\binom{\aleph_2}{\aleph_1} \rightarrow \binom{\aleph_0}{\aleph_0}^{1,1}_{\aleph_0}$.
\end{theorem}

\par\noindent\emph{Proof}. \newline
Let $\mathcal{J}$ be a uniform $\sigma$-complete $(\aleph_2,\aleph_1,\aleph_0)$-saturated ideal over $\aleph_1$.
Let $c:\aleph_2\times\aleph_1\rightarrow\aleph_0$ be any coloring.
For each $\alpha\in\omega_2$ we choose $n(\alpha)\in\omega$ so that $S_\alpha = \{\beta\in\omega_1: c(\alpha,\beta)=n(\alpha)\}\in \mathcal{J}^+$. This is possible since $\omega_1\in\mathcal{J}^+$ and $\mathcal{J}$ is $\sigma$-complete. Without loss of generality, there is a fixed color $m\in\omega$ such that $n(\alpha)=m$ for every $\alpha\in\omega_2$.

For every $\nu\in\omega_2$ let $\mathcal{S}_\nu = \{S_\alpha: \nu\leq\alpha\in\omega_2\}$. Let $\mathcal{T}_\nu = \{S_{\alpha^\nu_\tau}:\tau\in\omega_1\}$ be a subset of $\mathcal{S}_\nu$ such that $\bigcap_{\tau\in y}S_{\alpha^\nu_\tau}\in \mathcal{J}^+$ for every $y\in[\omega_1]^\omega$.
We shall use only one collection of the form $\mathcal{T}_\nu$, but we would like to choose it sufficiently high in the ordinals of $\omega_2$.
Let $\mathcal{W} = \{\mathcal{T}_\nu:\nu\in\omega_2\}$.

Let $\langle M_\eta:\eta\leq\omega\rangle$ be an increasing and continuous sequence of elementary submodels of $\mathcal{H}(\chi)$ for an appropriate regular cardinal $\chi$, such that $|M_\eta|=\aleph_1, \omega_1+1\subseteq M_\eta$ for every $\eta\leq\omega$. We may assume that $c,\mathcal{J},\mathcal{S}, \mathcal{W}\in M_0$ and that $\zeta<\eta\leq\omega \Rightarrow M_\zeta\in M_\eta$. We emphasize that the closure of the $M_\eta$-s under $\omega$-sequences is not assumed, as we wish to incorporate the cases of $2^\omega>\omega_1$ in our theorem. Let $\delta_0 = \sup(M_\omega\cap\omega_2)$.

We fix a collection $\mathcal{T}_\nu$ so that there is an element $S_{\alpha_\tau}\in\mathcal{T}_\nu$ for which $\alpha_\tau\geq\delta_0$. This is possible as we can step arbitrarily high with the starting point of $\mathcal{T}_\nu$. We denote the fixed collection by $\mathcal{T}$ and we enumerate its elements by $\{T_\alpha:\alpha<\omega_1\}$. Let $\delta$ be the ordinal for which $S_{\alpha_\tau}=T_\delta$.
We shall define two $\omega$-sequences of ordinals, namely $\langle\alpha_n:n\in\omega\rangle$ and $\langle\beta_n:n\in\omega\rangle$. The $\beta_n$-s will be ordinals in $\omega_1$, while the $\alpha_n$-s in $\omega_2$. They will form the monochromatic product which we try to create. The definition is done by induction on $n\in\omega$.

At the stage of $n=0$ we let $\beta_0 = \min(T_\delta)$.
By elementarity, there is some $\zeta\in M_1- M_0$ such that $\beta_0\in T_\zeta$, so we choose such an ordinal as $\alpha_0$. Arriving at $n>0$, let $\beta_n = \min(\bigcap_{\ell<n}T_{\alpha_\ell}\cap T_\delta- \{\beta_\ell:\ell<n\})$. Notice that $\beta_n$ is well-defined, by saturation.
Let $B_n = \{\beta_\ell:\ell\leq n\}$, so $B_n\in M_\eta$ for each $\eta\leq\omega$. Observe that $B_n\subseteq T_\delta$, so by elementarity there is an ordinal $\zeta\in M_{n+1}-M_n$ such that $B_n\in T_\zeta$. We choose such an ordinal to be $\alpha_n$.

Define $A = \{\alpha_n:n\in\omega\}, B = \{\beta_n:n\in\omega\}$. By the construction, $\beta_\ell\in T_{\alpha_n}$ for every $\ell,n\in\omega$. Hence $c''(A\times B) = \{m\}$ and the proof is accomplished.

\hfill \qedref{thmweaklaver}

The above theorems yield the following:

\begin{corollary}
\label{coranswer} Assume Martin's axiom and $2^\omega>\aleph_2$. Assume, further, that Chang's conjecture fails.
Then there are no weakly Laver ideals over $\aleph_1$.
\end{corollary}

\hfill \qedref{coranswer}

\begin{remark}
\label{rem0} Several comments.
\begin{enumerate}
\item [$(\alpha)$] The above proof generalizes to any $\kappa$ so that $\kappa = \kappa^{<\kappa}$. Namely, if there is a weakly Laver ideal over $\kappa^+$ then $\binom{\kappa^{++}}{\kappa^+} \rightarrow \binom{\kappa}{\kappa}^{1,1}_\kappa$.
\item [$(\beta)$] The saturation degree needed for Theorem \ref{thmweaklaver} is just $(\aleph_2,\aleph_1,<\omega)$-saturation. Hence the above corollary can be strengthened.
\item [$(\gamma)$] It seems that the proof of Corollary \ref{cortodorcevic} cannot be generalized, as it is, to higher versions of Martin's axiom.
The reason is the typical problematic issue which arises in this context. That is, higher versions of Martin's axiom require a strong form of the chain condition and a strong form of completeness.
\end{enumerate}
\end{remark}

We conclude with a couple of open problems. The first question is motivated by the last part of the above remark.

\begin{question}
\label{q0} Is it consistent that Baumgartner's axiom holds, $2^{\aleph_1}\geq\aleph_3$ and there exists a weakly Laver ideal over $\aleph_2$?
\end{question}

We used Martin's axiom in order to exclude the existence of weakly Laver ideals over $\aleph_1$. However, we had to assume that $2^\omega>\aleph_2$.
The reason is that the negative relation obtained by Todor\v{c}evi\'c requires $\aleph_2$-many dense sets in order to force its existence.

\begin{question}
\label{q1} Is it consistent that Chang's conjecture fails, Martin's axiom holds, $2^\omega = \omega_2$ and there exists a weakly Laver ideal over $\aleph_1$?
\end{question}

\newpage

\section{Acknowledgements}

I am grateful to the anonymous referee for the excellent work on my paper.
I also thank Yair Hayut for a very helpful discussion.

\newpage

\bibliographystyle{alpha}
\bibliography{arlist}

\end{document}